\documentclass[Svgc]{Svgc}

\usepackage{theorem}
\usepackage{graphics}
\usepackage{amsfonts}
\usepackage{latexsym}
\usepackage{epsfig}
\journalname{Graphs and Combinatorics}

\spnewtheorem{defi}{Definition}{\bf}{\rm}
\spnewtheorem{exam}{Example}{\bf}{\rm}

\begin{document}

\title{The Link Component Number of Suspended Trees}
\author{Toshiki Endo}

\institute{Jiyu Gakuen College, 1-8-15 Gakuen-cho, Higashikurume-shi, Tokyo 203-8521, Japan\\
e-mail: end@prf.jiyu.ac.jp}

\maketitle 

\begin{abstract}
This paper provides a relationship between a geometric structure of a suspended tree and the number of link components of the associated link diagram. 
\end{abstract}

\begin{keyword}
knot, link, link component number, suspended tree, planar graph
\end{keyword}

\section{Research Motivation and Early Studies}
We use standard terminology and notation of knot theory and graph theory, see for example \cite{A} and \cite{D}, respectively. Graphs considered in this paper are assumed to be embedded in the 2-sphere ${\mathbb S^2}$, that is, all graphs are {\it plane} graphs. Plane graphs are often used as a research tool in knot theory. This is because there is a one-to-one correspondence between a link diagram and an edge-signed plane graph. 

\begin{figure}[http]
\includegraphics[width=50mm]{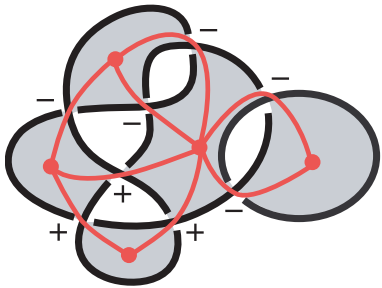}
\caption{A link diagram $L$ and its graph $G_L$}
\label{link_and_graph}
\end{figure}

Let $L$ be a link diagram in the 2-sphere ${\mathbb S^2}$. 
Suppose first that $L$ has at least one crossing. Regarding $L$ as a 4-regular plane graph, we color the faces black and white. From this coloring, we get an edge-signed plane graph $G_L$, where its vertices are the black faces and two vertices are joined by an edge if they share a crossing of $L$. Each edge is given a plus or minus sign according to the over/under information of the crossing. See Fig. \ref{link_and_graph}. 
If $L$ has no crossings, then $G_L$ has no edges, and the number of vertices is equal to the number of link components. 
Conversely, from the edge-signed plane graph $G_L$, we can restore the link diagram $L$ by considering the medial graph of $G_L$. 

Our research interest is to explore when a graph represent a knot and when a link. In general, how can we determine the number of components of a link diagram from the associated plane graph? For this purpose, we may only consider the underlying unsigned graph, since ignoring signs of edges does not change the number of components. 

\begin{defi}
\label{def:linkcompnum}
For a plane graph $G$, the number of components of the associated link diagram is called the {\it link component number} of $G$, and is denote by $l(G)$. 
\end{defi}


\begin{problem}
\label{prob}
For a plane graph $G$, find a method to determine the link component number $l(G)$. 
\end{problem}

\begin{exam}
\label{exam:tree}
Let $G$ be a tree, then $l(G)$ is equal to one. The proof is as follows. 
Recall that each degree one vertex of a tree is called a {\it leaf}, and any tree other than $K_1$ has at least two leaves. 
At each leaf, we untwist the string of the link diagram as the Reidemeister move I and contract the incident edge, and we get a smaller tree with the same link component number. By repeating this operation, any tree can reach $K_1$, and we can see that $l(G)$ is equal to one. 
\hfill$\Box$
\end{exam}

There are several early studies along this line although they possibly have slightly different expressions. 
See \cite{EHK} \cite{J} \cite{M} \cite{NW} for examples. 
The most noteworthy result is the following, and this may be a solution of the problem above. 

\begin{theorem}
\emph{(Schwarzler-Welsh \cite{SW})}
\label{thm:tutte}
Let $T(G,x,y)$ be the Tutte polynomial of a plane graph $G$. Then it holds that $\displaystyle{T(G,-1,-1)=}$ $\displaystyle{(-1)^{|E(G)|}(-2)^{l(G)-1}}$.
\hfill$\Box$
\end{theorem}

Our aim is now to find out a relation between geometric structures and the link component number of a plane graph. 

\section{Suspended Trees and their Link Component Numbers}
In this paper, we consider a certain extension of a tree, and completely determine its link component number. In what follows, we may assume that a tree has at least one edge, thus it is not $K_1$. Recall by Example \ref{exam:tree} that the link component number of a tree is equal to one, and this arises from the existence of a leaf. 

\begin{defi}
\label{def:suspetree}
The graph generated by a tree $T$ by adding a new vertex $v$ and new edges joining $v$ and all of the leaves of $T$ is called the {\it suspended tree}, and is denoted by $S_T$. 
\end{defi}

We begin with the following observation. 

\begin{proposition}
\label{prop:hasu}
For a tree $T$, the link component number of the suspended tree $S_T$ is at most the number of the leaves of $T$.
\end{proposition}

\begin{proof}
Note that $S_T-v$ has only one string. 
Consider the surrounding of the new vertex $v$. 
We may add up at most the number of arcs appeared here to the link component number of $S_T$, and this number equals to the number of the leaves of $T$. Thus, the proposition follows. 
\hfill$\Box$
\end{proof}

\begin{exam}
\label{exam:star}
Let $G$ be a star graph $K_{1,n}$, then $l(G)$ is equal to the number of the leaves $n$. In this case, it holds that $S_{T}=K_{2,n}$ and in the proof of Proposition \ref{prop:hasu}, the number of the arcs appeared in the surrounding of $v$ is equal to $n$ and the arcs make mutually different link components. 
\hfill$\Box$
\end{exam}


\section{The invariance under embeddings}
In this section, we discuss about the invariance of the link component number under planar embeddings. 

Let $G$ be a planar graph, and let $f_1$ and $f_2$ be two embeddings of $G$ into ${\mathbb S^2}$. Then plane graphs $f_1(G)$ and $f_2(G)$ are {\it equivalent} if there exists a homeomorphism of ${\mathbb S}^2$ onto itself which maps $f_1(G)$ into $f_2(G)$. From the well-known Whitney's theorem \cite{W}, if $G$ $3$-connected, then $f_1(G)$ and $f_2(G)$ are equivalent. 

\begin{figure}[htbp]
\includegraphics*[width=150mm]{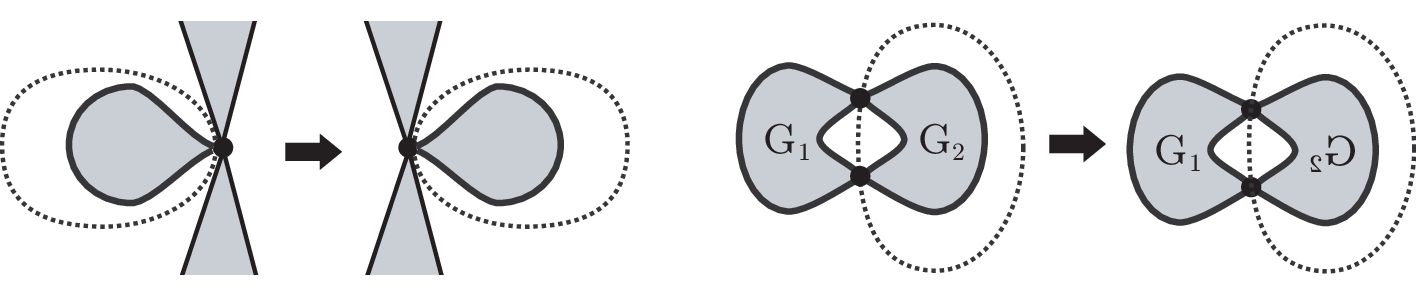}
\caption{a local jump and a local reversion}
\label{flipping}
\end{figure}

For the case that $G$ is not $3$-connected, Negami clarified the difference of the embeddings. In order to describe the statement, we review the definitions. 
Let $f$ be an embedding of a planar graph $G$ into ${\mathbb S^2}$. A {\it local jump} is defined as follows. Let ${\mathbb D}_1$ and ${\mathbb D}_2$ be two disks in ${\mathbb S^2}$ such that $\partial{\mathbb D}_1\cap f(G)=\partial{\mathbb D}_2\cap f(G)$ is a vertex $f(v)$. Then a local jump is the modification of $f$ into $\tau\circ f$ where $\tau$ is an orientation-preserving homeomorphism of ${\mathbb D}_1$ into ${\mathbb D}_2$. 
A {\it local reversion} is defined as follows. Let ${\mathbb D}$ be a disk in ${\mathbb S^2}$ and let $\tau$ is an orientation-reversing homeomorphism of ${\mathbb D}$ onto itself such that each point $x\in\partial{\mathbb D}\cap f(G)$ is a vertex of $f(G)$ and $\tau(x)=x$. We may assume that $\partial{\mathbb D}\cap f(G)$ consists of at most two vertices. Then a local reversion is the alternation of $f$ into $\tau\circ f$. See Fig. \ref{flipping}. 

\begin{theorem}
\emph{(Negami \cite{N})}
\label{thm:negami}
Let $G$ be a planar graph, and let $f_1$ and $f_2$ be two embeddings of $G$ into ${\mathbb S^2}$. Then $f_1(G)$ and $f_2(G)$ may be equivalent or $f_1(G)$ can be transformed into $f_2(G)$ by a finite sequence of local jumps and local reversions. 
\hfill$\Box$
\end{theorem}

Now we prove the following. 

\begin{theorem}
\label{thm:umekomi}
Let $G$ be a planar graph, and $f_1$ and $f_2$ be two embeddings of $G$ into ${\mathbb S^2}$. Then the link component numbers of $f_1(G)$ and $f_2(G)$ are equal. 
\end{theorem}

\begin{proof}
From Theorem \ref{thm:negami}, we may only show that the theorem is true for the case that $f_2(G)$ is obtained from $f_1(G)$ by performing a local jump and the case by a local reversion. 

(Case 1: local reversion)\quad
Suppose first that $\partial{\mathbb D}\cap f_1(G)$ consists of two vertices $f_1(u_1)$ and $f_1(u_2)$. Then the associated link diagram of $f_1(G)$ meets $\partial{\mathbb D}$ at four points. Let $t_1, t_2, t_3, t_4$ be the points in a counterclockwise direction on $\partial{\mathbb D}$. 
Let ${\mathbb D}\cap f_1(G)=f_1(G_2)$ and $({\mathbb S}^2-int{\mathbb D})\cap f_1(G)=f_1(G_1)$. 
Let $a_i$ (resp. $b_i$) be the string of the associated link diagram of $f_1(G_1)$ (resp. $f_1(G_2)$) which contains $t_i$ for $i=1,2,3,4$. As for the associated link diagram of $f_1(G_1)$, there are three possible combinations: (a1) $a_1=a_2$ and $a_3=a_4$, (a2) $a_1=a_3$ and $a_2=a_4$, and (a3) $a_1=a_4$ and $a_2=a_3$. And this is the case with the associated link diagram of $f_1(G_2)$, and we similarly denote them by (b1), (b2), and (b3). 

Although we have $a_i=b_i$ for $i=1,2,3,4$ in the associated link diagram of $f_1(G)$, we have $a_1=b_2, a_2=b_1, a_3=b_4, a_4=b_3$ in the associated link diagram of $f_2(G)$. Thus the combination of the strings changes by the local reversion, but the link component numbers does not change in each case. These are arranged in the table bellow. 

The subcases that $\partial{\mathbb D}\cap f_1(G)$ consists of one vertex and null can be similarly and more easily shown. 

\begin{table}[htbp]
\label{tab:string} 
\begin{tabular}{l|lll}
\hline
  & b1 & b2 & b3 \\
\hline
a1& 2 & 1 & 1 \\
a2& 1 & 2 & 1 \\
a3& 1 & 1 & 2 \\
\hline
\end{tabular}
\caption{Numbers of strings that meet $\partial{\mathbb D}$}
\end{table}

(Case 2: local jump)\quad
In this case, the associated link diagram of $f_i(G)$ meets $\partial{\mathbb D}_i$ at two points for $i=1, 2$. But the two points are contained in the same link component inside and outside ${\mathbb D}_i$, respectively. Thus the link component numbers does not change. 
\hfill$\Box$
\end{proof}

By Theorem \ref{thm:umekomi}, we may not distinguish between a planar graph and its plane graph as for the link component number. 
\bigskip

Next theorem is also noteworthy, and is used in Corollary \ref{cor:forest}. 

\begin{theorem}
\label{thm:1-sep}
Let $G$ be a planar graph and $u$ be a cut vertex such that $G=G_1\cup G_2$ and $\{u\}=G_1\cap G_2$. 
Then, the link component number of $G$ is equal to the sum of two link component numbers of $G_1$ and $G_2$ minus one. 
\end{theorem}

\begin{proof}
Let $L_1$ and $L_2$ be two link diagrams arising from $G_1$ and $G_2$, respectively. Then, we obtain a link diagram of $G$ by amalgamating one string of $L_1$ and one string of $L_2$ at $u$. Thus, the theorem follows. 
\hfill$\Box$
\end{proof}

\section{Key Lemmas and the Theorem}

First we observe the following. 

\begin{lemma}
\label{lem:contract}
Let $T$ be a tree, and suppose that $T$ has a vertex $v$ of degree two. Let $x, y$ be the neighbors of $v$, and suppose that neither $x$ nor $y$ is a leaf of $T$. Let $T'$ be the tree obtained from $T$ by contracting $vx$ and $vy$. 
Then the link component numbers of $S_T$ and $S_{T'}$ are equal. 
\end{lemma}

\begin{proof}
This is obvious if we see how strings of the corresponding link diagram behave. Imagine the Reidemeister move II. 
\hfill$\Box$
\end{proof}

\begin{lemma}
\label{lem:unavoidable}
Let $T$ be a tree which is not $K_2$ such that each vertex of degree two, if any, is adjacent to a leaf of $T$. Then, $T$ has at least one of the following vertex $w$: (Type I) $w$ is adjacent to a leaf and a degree two vertex which is adjacent to a leaf, (Type II) $w$ is adjacent to two degree two vertices each of which is adjacent to a leaf, and (Type III) $w$ is adjacent to two leaves. 
\end{lemma}

\begin{proof}
We show by contradiction. Suppose $T$ has no vertices of Type I, II, or III. Let $x, y, z$ be numbers of the vertices of degree one, two, and three or more, respectively. 

Then, for any leaf $l$, there exists a vertex $w$ with $\deg_Tw\geq 3$ such that $w$ is adjacent to $l$ or a degree two vertex which is adjacent to $l$. In this situation, we say that $v$ has $l$ as a {\it root}. Then, for any vertex $w$ with $\deg_Tw\geq 3$, $w$ has at most one leaf as a root. Thus, it holds that $z\geq x$. 

Since $T$ is a tree, the number of the edges is $x+y+z-1$, and $x+y+z-1\geq (x+2y+3z)/2$ by the handshaking lemma. Thus we have $x-2\geq z$, a contradiction. \hfill$\Box$
\end{proof}

\begin{lemma}
\label{lem:typeI}
Suppose that a tree $T$ has a vertex $w$ of Type I. 
Let $T'$ be the tree obtained from $T$ by deleting a length two path from $w$ to a leaf. Then, the link component numbers of $S_T$ and $S_{T'}$ are equal. 
\end{lemma}

\begin{figure}[htbp]
\includegraphics*[width=60mm]{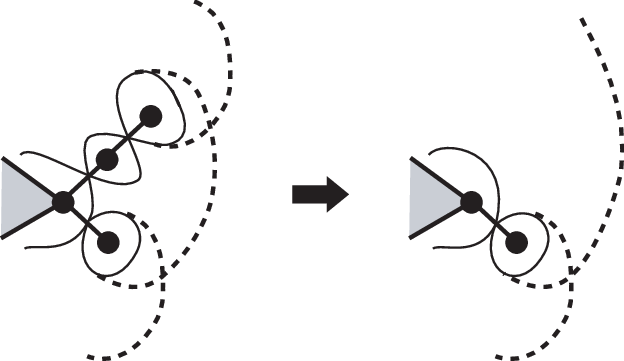}
\caption{Deletion of a path of length two which is incident to a vertex of Type I}
\label{delete2}
\end{figure}

\begin{proof}
If $\deg_Tw\geq 3$, then this is obvious from Fig. \ref{delete2}, where dotted lines denote strings around the new vertex $v$ of $S_T$. 
In the case that $\deg_Tw=2$, although non-leaf $w$ of $T$ changes a leaf of $T'$, it is not difficult to check the link component numbers of $S_T$ and $S_{T'}$ are equal. 
\hfill$\Box$
\end{proof}

\begin{lemma}
\label{lem:typeII}
Suppose that a tree $T$ has a vertex $w$ of Type II. 
Let $T'$ be the tree obtained from $T$ by deleting two length two paths each of which connects $w$ and a leaf, and if $\deg_Tw=3$, then we delete $w$ furthermore. Then, the link component numbers of $S_T$ and $S_{T'}$ are equal. 
\end{lemma}

\begin{figure}[htbp]
\includegraphics*[width=130mm]{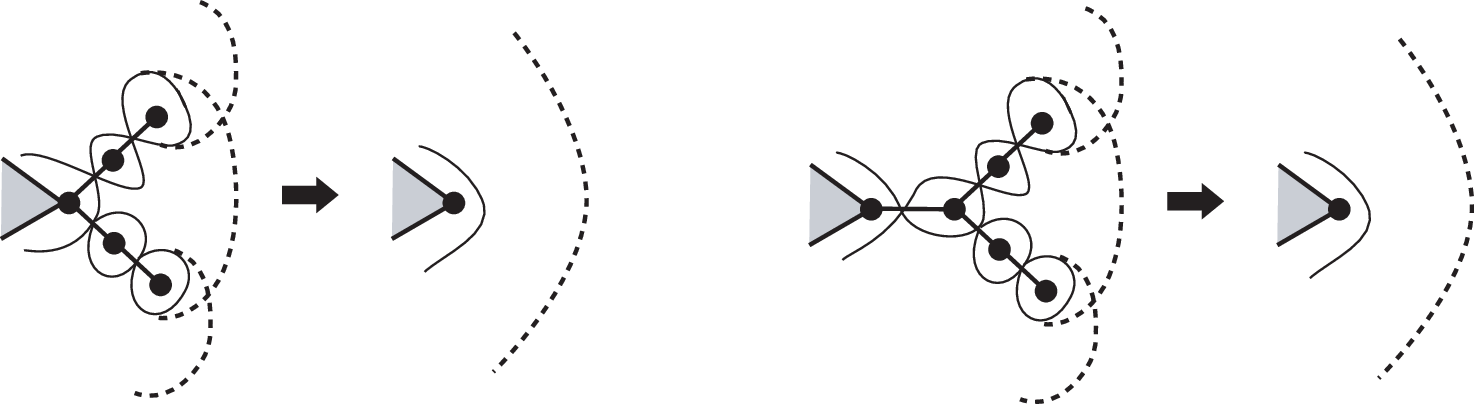}
\caption{Deletion of two paths of length two which are incident to a vertex of Type II}
\label{delete2w}
\end{figure}

\begin{proof}
This is obvious from Fig. \ref{delete2w}. 
In the case that $\deg_Tw=3$, there is a possibility that the third neighbor $x$ of $w$ has degree two in $T$ and consequently $x$ might have degree one in $T'$. But in that case, the neighbor of $x$ must be a leaf, and hence the lemma follows. 
\hfill$\Box$
\end{proof}

\begin{lemma}
\label{lem:typeIII}
Suppose that a tree $T$ has a vertex $w$ of Type III. 
Let $T'$ be the tree obtained from $T$ by deleting two leaves. 
Then, the link component number of $S_{T'}$ is equal to the link component number of $S_{T}$ minus one. 
\end{lemma}

\begin{figure}[htbp]
\includegraphics*[width=55mm]{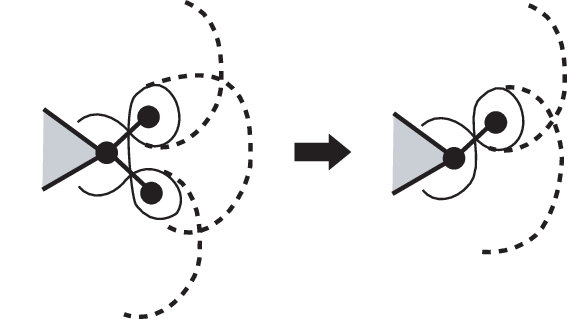}
\caption{Deletion of an edge which is incident to a vertex of Type III}
\label{matsuba}
\end{figure}

\begin{proof}
If $\deg_Tw\geq 3$, then this is obvious from Fig. \ref{matsuba}. 
In the case of $\deg_Tw=2$, although non-leaf $w$ of $T$ changes a leaf of $T'$, it is not difficult to check the lemma holds. 
\hfill$\Box$
\end{proof}

The Main Theorem in this paper is the following. 

\begin{theorem}
\label{thm:main}
Any tree can be transformed into $K_2$ by performing the operations as in Lemmas \ref{lem:contract}, \ref{lem:typeI}, \ref{lem:typeII} and \ref{lem:typeIII}. Moreover, the link component number of $S_T$ is equal to the number of times of the operations of Lemma \ref{lem:typeIII} plus one. 
\end{theorem}

\begin{proof}
We show the theorem by a mathematical induction on the number of the edges of $T$. The first step is the case that $T=K_2$, and it is trivial. The second step immediately follows from Lemmas \ref{lem:contract}--\ref{lem:typeIII}.
\hfill$\Box$
\end{proof}

Theorem \ref{thm:main} suggests that the link component number of a suspended tree depends on essentially the number of $K_{1,2}$'s appeared in the operations above, where the two vertices of one partite set are leaves. Comparing to the computation by using the Tutte polynomial (Theorem \ref{thm:tutte}), the computational effort may be relatively small. 
\bigskip

We conclude this section by mentioning two corollaries of our theorem. 

\begin{defi}
\label{def:suspefore}
The {\it suspended forest} is the graph obtained from a forest by adding a new vertex and new edges joining the vertex and the leaves of the forest. 
\end{defi}

\begin{corollary}
\label{cor:forest}
Let $F$ be a forest consisting of $n$ trees. Then, the link component number of the suspended forest is equal to the sum of the link component numbers of the constituent suspended trees plus $n$ minus one. 
\end{corollary}

\begin{proof}
This is an immediate consequence of Theorem \ref{thm:1-sep}. 
\hfill$\Box$
\end{proof}

\begin{corollary}
\label{cor:notknot}
Let $T$ be a tree other than $K_2$, and has no vertices of degree two. Then, the link component number of the suspended tree is two or more. \hfill$\Box$
\end{corollary}

\begin{proof}
This is because $T$ must have a vertex of type III. 
\hfill$\Box$
\end{proof}

Thus, in contrast to the fact that a tree always represents a knot, a suspended tree does not represent a knot if the tree part has no vertices of degree two. 

\section{Concluding Remark}
We shall briefly describe a relationship between suspended trees and knot theory. A link $L$ is called {\it arborescent} (or {\it algebraic}) if $L$ is formed by taking the numerator closure of an tangle obtained by additions and multiplications of rational tangles. It is known that a link $L$ is arborescent if and only if the associated graph $G_L$ has a vertex $v$ such that $T=G_L-v$ is a tree \cite{MP}. Thus, as an application of Theorem \ref{thm:main}, we may determine the link component number of an arborescent link via its graph. 
In fact, if $v$ is adjacent to some vertex of $T$, the number of edges between them may be assumed to be at most one, since we can delete parallel two edges by considering the Reidemeister move II. If $v$ is adjacent to a non-leaf vertex $x$ of $T$, then we may consider the tree obtained from $T$ by adding a length two path to $x$. 

\begin{acknowledgement}
The author is grateful to the referee for many helpful suggestions that lead to improvements. 
\end{acknowledgement}

\end{document}